\documentclass[10pt]{article}

\usepackage{amsmath}
\usepackage{amssymb}
\usepackage{indentfirst}
\usepackage{graphics} 
\usepackage{color}
\makeatletter

\@addtoreset{equation}{section}
\makeatother
\def\<{\langle}
\def\>{\rangle}

\newtheorem{lem}{Lemma}[section]
\newtheorem{theo}{Theorem}[section]
\newtheorem{rem}{Remark}[section]

\makeatletter
   
   \@addtoreset{equation}{section}
\makeatother

\begin{document}
\title{\bf A note on local energy decay results\\for wave equations with a potential}
\author{Ryo Ikehata\thanks{Corresponding authr: ikehatar@hiroshima-u.ac.jp} \\ {\small Department of Mathematics, Division of Educational Sciences}\\ {\small Graduate School of Humanities and Social Sciences} \\ {\small Hiroshima University} \\ {\small Higashi-Hiroshima 739-8524, Japan} %\\ and \\Ruikson S. O. Nunes\thanks{ruiksonsillas@hotmail.com} \\{\small UFMT-Federal University of Mato Grosso, ICET} \\{\small Department of Mathematics} \\ {\small 78060-900, Cuiab\'a, Brazil}
}
\date{}
\maketitle
%%%%%%%%%%%%%%%%%%%%%%%%%ހabstractހ%%%%%%%%%%%%%%%%%%%%%%
\begin{abstract}
In this paper we consider the local energy decay result for wave equations with a short-range potential. It is important to note that one never uses a finite speed of propagation property unlike the historical previous papers. The essential parts of analysis are in getting $L^{2}$-bound of the solution itself, and deriving the weighted energy estimates. In this paper we only use a simple multiplier method to treat the variable coefficient case, and do not rely on any resolvent estimates.  
\end{abstract}
%%%%%%%%%%%%%%%%%%ހintroductionހ%%%%%%%%%%%%%%%%%%%%%%%%
\section{Introduction and Statement of Results}
\footnote[0]{Keywords and Phrases: Wave equation; short-range potential; exterior mixed problem, non-compact support, initial data, local energy, algebraic decay.}
\footnote[0]{2010 Mathematics Subject Classification. Primary 35L05; Secondary 35B40, 35L30.}

In this paper, we are concerned with the following Cauchy problem with a potential:

\begin{equation}\label{01}
u_{tt}- \Delta u + V(x)u = 0, \quad t > 0,\quad x \in \Omega,
\end{equation}
\begin{equation}\label{02}
u(0,x)=u_{0}(x),\hspace{0,50cm}u_{t}(0,x) = u_{1}(x)\quad x \in \Omega,
\end{equation} 
\begin{equation}\label{03}
u(t,x) = 0\quad t > 0, \quad x \in \partial\Omega,
\end{equation} 
where $\Omega \subset {\bf R}^{n}$ is an exterior domain with smooth compact boundary $\partial\Omega$, and one assumes $0 \notin \bar{\Omega}$. Furthermore, let $\rho_{0} > 0$ be a real number such that
$\partial\Omega \subset B_{\rho_{0}}$, and one assumes that\\
{\bf (A-1)} the obstacle ${\cal O} := {\bf R}^{n}\setminus \bar{\Omega}$ is star-shaped relative to the origin, that is, $x\cdot \nu(x) \leq 0, \; x \in \partial\Omega$ where $\nu(x)$ is the unit exterior normal at the point $x \in \partial \Omega.$

The potential function $V \in BC^{1}(\bar{\Omega})$ satisfies $V (x) \geq 0$ ($x \in \Omega$) and\\ 
{\bf (A-2)}\,\,$\displaystyle{\frac{1}{2}}\left(x\cdot\nabla V(x)\right) + V(x) \leq 0$ for all $x \in \bar{\Omega}$.\\
Note that functions and solutions treated in this paper are all real-valued.\\
\noindent
{\bf Example 1.} One can present a typical example for $V(x)$ satisfying the assumption {\bf (A-2)} as follows:
\[V(x) = V_{0}\vert x\vert^{-\alpha},\quad V_{0} > 0,\]
where $\alpha \geq 2$. This example shows that $V(x)$ is a short-range potential.

\begin{rem}{\rm It seems quite important for this type of problems whether the case $\alpha = 2$ can be included or not as an example of $V(x) = V_{0}\vert x\vert^{-\alpha}$ (cf., \cite{GV} and \cite{V}).}
\end{rem}
\begin{rem}{\rm In the case of radial function $V(x) = V(r)$ for $r := \vert x\vert$ the condition {\bf (A-2)} can be replaced by
\[rV'(r) \leq -2V(r),\quad x \in \bar{\Omega}.\]
}
\end{rem}
{\bf Example 2.}\,If ${\cal O} := B_{2}$, then we can choose $V(x) := V_{0}e^{-\vert x\vert}$ with $V_{0} > 0$ as the second example.\\

\begin{rem}{\rm The condition (8) assumed on the paper by \cite{C-Ruy} is closely related to (A-2):
\begin{equation}\label{ike-501}
\frac{1}{2}\left(x\cdot\nabla V(x)\right) + V(x) \leq \frac{\gamma}{2}V(x),\quad x \in \bar{\Omega}
\end{equation}
with $\gamma \in [0,1)$, that is, if we choose $\gamma = 0$ formally, the condition \eqref{ike-501} of \cite{C-Ruy} is (A-2). So in this sense, the assumption \eqref{ike-501} in \cite{C-Ruy} is weaker than (A-2). However the local energy decay in \cite{C-Ruy} was obtained on the strong assumption that  $q(x)$has  compact support and the finite speed of propagation. The work in \cite{C-Ruy} is to the system of elastic waves in $\mathbb{R}^3$ but of course the results hold to the wave equation with the same type of potential.}
\end{rem}

One defines the total energy for the equation \eqref{01} by
\begin{equation}\label{E1}
E(t) := \int_{\Omega} e(t,x)dx := \frac{1}{2}\int_{\Omega}\left(\vert u_{t}(t,x)\vert^{2} + \vert\nabla u(t,x)\vert^{2} + V(x)\vert u(t,x)\vert^{2} \right)dx.
\end{equation}

Then, under the assumptions {\bf (A-1)} and {\bf (A-2)} it is known that for each initial data $[u_{0},u_{1}] \in H_{0}^{1}(\Omega) \times L^{2}(\Omega)$ the problem \eqref{01}-\eqref{03} has a unique weak solution $u \in C([0,\infty); H_{0}^{1}(\Omega)) \cap C^{1}([0,\infty); L^{2}(\Omega))$ satisfying the energy identity 
$$E(t) = E(0).$$

Our main purpose is to study the local energy decay problem of the equation \eqref{01} with the short-range potential $V(x)$. Here, for each $R > 0$ the local energy $E_{R}(t)$ can be defined as follows:
\[E_{R}(t) := \frac{1}{2}\int_{B_{R}\cap \Omega}\left(\vert u_{t}(t,x)\vert^{2} + \vert\nabla u(t,x)\vert^{2} + V(x)\vert u(t,x)\vert^{2} \right)dx.\]

First of all, from the multiplier method point of views one can cite the celebrated work due to C. Morawetz \cite{maraw}. There the local energy decay result can be studied by basing on the so-called Morawetz identity to the equation \eqref{01} with $V(x) \equiv 0$. In fact, Morawetz derived the fact that $E_{R}(t) = O(t^{-1})$ ($t \to \infty$) under the stronger geometrical condition on the obstacle (star-shaped obstacle case). To obtain such result, the author in \cite{maraw} assumed that the initial data have compact support to use the finite speed of propagation property. One of the essential ingredient is to get the $L^{2}$-bound of the solution itself, and the estimate for the solution of the corresponding Poisson equation played a crucial role in \cite{maraw}. In particular, in the three dimensional case it seems easy to prove that the defined local energy, by using Huygens principle, decays exponentially fast (see \cite{morast}), and the non-trapping obstacle can be treated in the study as a generalization on the obstacle. \\

After Morawetz, several efforts to remove the compactness assumptions on the support of the initial data are devoted by \cite{Zac-2}, \cite{Mur} and \cite{Ikehat, Ikehat2, Ikken}. The authors in \cite{Zac-2}, \cite{Mur} and \cite{Ikehat, Ikehat2, Ikken} have employed the multiplier method, and the decay rate $O(t^{-2})$ and the integrability of the local energy with respect to the time variable $t$ can be derived by \cite{Zac-2}, \cite{Mur} under the rather stronger weight condition on the initial data, while the decay rate $O(t^{-1})$ of the local energy has been got under weaker weight condition on the initial data due to \cite{Ikehat, Ikehat2, Ikken} (see also \cite{CI} for the variable coefficient case with Lipschitz wave speeds). It should be mentioned that the weight condition (as $\vert x\vert \to \infty$) assumed on the initial data seems to be the weakest within the known results. One can also cite the other related deep investigations on the topics of the local energy decay due to \cite{B}, \cite{LP}, \cite{R}, \cite{S}, \cite{ST}, and \cite{Va} under the conditions of the compact support initial data. In particular, in \cite{AGPP} one dimensional wave equations with the variable coefficients are taken up to capture the exponential decay of the local energy. To the best of authors'knowledge, this research of \cite{AGPP} seems first to take up deeply the one dimension case. 
   
On the other hand, as for the same topics to the equation \eqref{01} a few results are known, and in particular, the author in \cite{V} has investigated the sharp local energy decay rates in the short-range case such that $V(x) = V_{0}<x>^{-\alpha}$ satisfying $\alpha > 2$. In fact, the author studies the Cauchy problem of \eqref{01} in ${\bf R}^{n}$ ($n \geq 3$), and get the decay rate $O(t^{-2})$. In that paper \cite{V}, the compactness of the support of the initial data is essentially used. In this connection, uniform weighted resolvent estimates are effectively employed in \cite{V}. Therefore, it seems unknown to consider the local energy decay problem to the equation \eqref{01} in the case when the compact support assumption is not assumed on the initial data. We develop our theory by using the multiplier method based on the expanded Morawetz identity. As a side note one should two references  \cite{C-Ruy} and \cite{T}. There the local energy decay problems are investigated for the elastic waves with time-independent potentials and wave equations with time-dependent potentials, respectively, however, the authors in  \cite{C-Ruy} and \cite{T} studied such problems in the framework of compact support assumptions to both potentials and initial data. It should be noted that the condition (A-2) in our paper is stronger
than that assumed in (8) of \cite{C-Ruy} (See Remark 1.2 ).

For later use, one defines a weight function $d_{n}(x)$ by
\[ d_{n}(x) = \left\{
    \begin{array}{ll}
    \displaystyle{\vert x\vert}&
    \qquad n \geq 3, \\[0.2cm]
    \displaystyle{\log(B\vert x\vert)}& \qquad n = 2,
    \end{array} \right.\]
with some constant $B > 0$ satisfying $B\inf\{\vert x\vert\,:\, x \in \Omega\} \geq 2$. 

The main result of this paper reads as follows.

\begin{theo}\label{t01} 
Let $\rho_0>0$ such that $\partial\Omega\subset{B}_{\rho_0}$.
Let $n\geq 2$ and assume {\bf (A-1)} and {\bf (A-2)}. Let $R > \rho_{0}$ be an arbitrary fixed number. If $(u_{0}, u_{1})\in C^{\infty}_{0}(\Omega)\times C_{0}^{\infty}(\Omega)$, then the unique smooth solution $u(t,x)$ to problem \eqref{01}-\eqref{03} satisfies
\[E_{R}(t) \leq \frac{CK_{0}}{t-R}, \quad t > R,\]
with some generous constant $C > 0$, where
\[K_{0} := (u_{1},u_{0}) + (u_{1},x\cdot\nabla u_{0}) + \sqrt{E(0)}(\Vert u_{0}\Vert + \Vert d_{n}(\cdot)u_{1}\Vert)\]
\[+ \int_{\Omega}(1+\vert x\vert)\left(\vert u_{1}(x)\vert^{2} + \vert\nabla u_{0}(x)\vert^{2} + V(x)\vert u_{0}(x)\vert^{2}  \right)dx.\]
\end{theo}

\begin{rem}{\rm In some sense the assumption (A-2) on the potential $V(x)$ is a technical condition, however it includes an important example $\alpha = 2$, that is, $V(x) = \vert x\vert^{-2}$ as a critical potential. An important fact is that such
singular potential is a unique type of potential that the perturbed wave equation retains the property of Huygen's principle in dimension $n=3$. For the perturbed wave equation with a regular potential, the Huygens' principle never holds (see \cite{C-Perla}). In that case, if the local energy decays then it decays exponentially in case of the validity of such a principle. }
\end{rem}

It should be emphasized that the constant $C > 0$ determined in Theorem \ref{t01} never depends on any $R > \rho_{0}$, and $R > \rho_{0}$ is independent of the size of support of the initial data. This implies that one never relies on the finite speed of propagation property as is usually discussed (cf., \cite{maraw}). This is our essential contribution, and the condition $(u_{0}, u_{1}) \in C^{\infty}_{0}(\Omega)\times C_{0}^{\infty}(\Omega)$ assumed on the initial data is not essential. By density one can discuss the same local energy decay in the framework of $H_{0}^{1}(\Omega)\times L^{2}(\Omega)$. Indeed, one can obtain the following refinement of Theorem \ref{t01}. We state it without proof.
\begin{theo}\label{t02} 
Let $n\geq 2$ and assume {\bf (A-1)} and {\bf (A-2)}. Let $R > \rho_{0}$ be an arbitrary fixed number. If $(u_{0}, u_{1})\in H_{0}^{1}(\Omega)\times L^{2}(\Omega)$, then the unique weak solution $u \in C([0,\infty); H_{0}^{1}(\Omega)) \cap C^{1}([0,\infty); L^{2}(\Omega))$ to problem \eqref{01}-\eqref{03} satisfies
\[E_{R}(t) \leq \frac{CK_{0}}{t-R}, \quad t > R,\]
with some generous constant $C > 0$, and $K_{0}$ defined in Theorem {\rm \ref{t01}}, provided that 
\[\Vert d_{n}(\cdot)u_{1}\Vert < +\infty, \quad \int_{\Omega}\vert x\vert\left(\vert u_{1}(x)\vert^{2} + \vert\nabla u_{0}(x)\vert^{2} + \vert u_{0}(x)\vert^{2}\right)dx < +\infty.\]
\end{theo}
\begin{rem}{\rm Unfortunately, the constant coefficient case $V(x) = m^{2}$ ($m > 0$) can not be included as an example. This is the so-called Klein-Gordon equation case, which seems more difficult to be treated by our method. The assumption {\rm (A-2)} may express a small perturbation from the pure wave equation case with $V(x) \equiv 0$. For the sharp local energy decay of the Klein-Gordon equation by using compactness assumptions on the initial data, one can cite the recent paper \cite{Nu}.}
\end{rem}  
\begin{rem}{\rm In the assumptions on the initial velocity $u_{1}(x)$ of Theorem 1.2, it is easy to see that in the case when $n = 2$ the condition $\Vert d_{2}(\cdot)u_{1}\Vert < +\infty$ can be absorbed into $\displaystyle{\int_{\Omega}}\vert x\vert\vert u_{1}(x)\vert^{2}dx < \infty$, while in the case of $n \geq 3$, $\Vert d_{n}(\cdot)u_{1}\Vert < +\infty$ implies $\displaystyle{\int_{\Omega}}\vert x\vert\vert u_{1}(x)\vert^{2}dx < \infty$.}
\end{rem}
Note that the concrete case $V(x) := V_{0} \vert x\vert^{-\alpha}$ with $\alpha \geq 2$ can be included as an example, and in this case from Theorem \ref{t01} one has
\[\frac{V_{0}}{2R^{\alpha}}\int_{B_{R}\cap \Omega}\vert u(t,x)\vert^{2}dx \leq \frac{1}{2}\int_{B_{R}\cap \Omega}V(x)\vert u(t,x)\vert^{2}dx \leq E_{R}(t) \leq \frac{CK_{0}}{t-R}, \quad t > R,\]  
so that one has also local $L^{2}$-decay result:
\begin{equation}\label{100}
\int_{B_{R}\cap \Omega}\vert u(t,x)\vert^{2}dx \leq \frac{2R^{\alpha}}{V_{0}}\frac{CK_{0}}{t-R}, \quad t > R.
\end{equation}

The decay result \eqref{100} is closely related to that of \cite[Theorem 1.2]{V}. In \cite{V} the critical case $\alpha =2$ cannot be included as an example.\\

{\bf Notation.}\, We denote the $L^{2}$-norm of $u \in L^{2}(\Omega)$ by $\Vert u\Vert$. We set $B_{R} := \{x \in {\bf R}^{n}\,:\,\vert x\vert < R\}$, and $(f,g) := \displaystyle{\int_{\Omega}f(x)g(x)dx}$ denotes the usual $L^{2}$-inner product of $f,g \in L^{2}(\Omega)$. We denote $f \in BC(\bar{\Omega})$ $\Leftrightarrow$ $f(x)$ is bounded and continuous in $\bar{\Omega}$. $f \in BC^{1}(\bar{\Omega})$ $\Leftrightarrow$ $f,\,\partial f/\partial x_{j} \in BC(\bar{\Omega})$ for $j = 1,2,\cdots, n$.\\

The rest of this paper is organized into three sections. The Section 2 is dedicated to expose some preliminary results which are used in the proof of Theorem \ref{t01}. In Section 3 one proves our main Theorem \ref{t01}.

%%%%%%%%%%%%%%%%%%%%%%%%%%%%%%%%%%%%%%%%%%%%%%%%%%%%%%%%%%%%%%%%%%%%%%%%%%%%%%%%%%%%%%%%%%%%%%%%%%%%%%%%%%%%%%%%%%%%%%%%%%%%%%%%%%%%%%%%%%%
\section{Preliminaries}

The following lemma is a kind of Morawetz identity for the equation \eqref{01}-\eqref{02} obtained by using the multiplier $m(u)=u_t+x\cdot\nabla u + \frac{n-1}{2}$. The Morawetz identity is the useful framework to have identities when one needs to obtain estimates on the solutions, at least to hyperbolic equations. In \cite{EVR} (Lemma 3.3) the authors obtaines identities with simple generalized multipliers of Morawetz type  to study the stabilization of solutions to the system of elastic waves with localized nonlinear dissipation.

\begin{lem}\label{10}\, Let $n \geq 2$, and $[u_{0},u_{1}] \in C_{0}^{\infty}(\Omega) \times C_{0}^{\infty}(\Omega)$. Then, the corresponding smooth solution $u(t,x)$ to problem \eqref{01}-\eqref{03} satisfies the following identity: for $t \geq 0$ it holds that
\[t E(t) + \frac{n-1}{2}\int_{\Omega}u_{t}(t,x)u(t,x)dx + \int_{\Omega}u_{t}(t,x)(x\cdot\nabla u(t,x))dx\]
\[-\int_{0}^{t}\int_{\Omega}\left(\frac{1}{2}(x\cdot\nabla V(x)) + V(x)\right)\vert u(s,x)\vert^{2}dx ds = J_{0} + \frac{1}{2}\int_{0}^{t}\int_{\partial\Omega}\left( \frac{\partial u}{\partial\nu} \right)^{2}\sigma\cdot\nu(\sigma)dS_{\sigma}ds,\]
where
\[J_{0} := \frac{n-1}{2}\int_{\Omega}u_{1}(x) u_{0}(x)dx + \int_{\Omega}u_{1}(x)(x\cdot \nabla u_{0}(x))dx,\]
and $\nu(\sigma)$ is the unit outward normal vector at each $\sigma \in \partial\Omega$.
\end{lem}

Furthermore, one needs the weighted energy estimate below. This is the modified version of that introduced originally by Todorova-Yordanov \cite{TY}(see also the Appendix in \cite{Ikken}).  For this we prepare the following notation for the pointwise total energy and the weight function, respectively.
\[e(t,x) :=  \frac{1}{2}\left(\vert u_{t}(t,x)\vert^{2} + \vert\nabla u(t,x)\vert^{2} + V(x)\vert u(t,x)\vert^{2}\right), \quad t > 0, \quad x \in \Omega,\]
and $\psi \in C^{1}([0,\infty)\times \bar{\Omega})$ satisfying $\psi_{t}(t,x) \ne 0$ for all $(t,x) \in [0,\infty)\times \bar{\Omega}$.
\begin{lem}\label{11}\, Let $n \geq 2$, and $[u_{0},u_{1}] \in C_{0}^{\infty}(\Omega) \times C_{0}^{\infty}(\Omega)$. Then, the corresponding smooth solution $u(t,x)$ to problem \eqref{01}-\eqref{03} satisfies the following identity:
\[0 = \frac{\partial}{\partial t}\left(\psi(t,x)e(t,x)\right) - \nabla\cdot\left(\psi(t,x)u_{t}(t,x)\nabla u(t,x) \right) - \frac{V(x)}{2}\vert u(t,x)\vert^{2}\psi_{t}(t,x)\]
\[-\frac{1}{2\psi_{t}(t,x)}\left\vert \psi_{t}(t,x)\nabla u(t,x)-u_{t}(t,x)\nabla\psi(t,x)\right\vert^{2} + \frac{\vert u_{t}(t,x)\vert^{2}}{2\psi_{t}(t,x)}\left(\vert\nabla \psi(t,x)\vert^{2}-\psi_{t}(t,x)^{2}\right),\,\,t > 0,\,\,x \in \Omega.\]
\end{lem}

To prove the  following $L^{2}$-estimate of the solution one may derive itself  by a similar method introduced in \cite{IM} (see also \cite[Lemma 2.2]{Ikken}) even for the equation \eqref{01} with potential $V(x)$. Since the proof relies on the Hardy inequality in the exterior domains for $n \geq 2$, the weight function $d_{n}(x)$ appears in the statement (see \cite{DS}). 
\begin{lem}\label{12}\, Let $n \geq 2$, and $[u_{0},u_{1}] \in C_{0}^{\infty}(\Omega) \times C_{0}^{\infty}(\Omega)$. Then, the corresponding smooth solution $u(t,x)$ to problem \eqref{01}-\eqref{03} satisfies the following estimate:
\[\Vert u(t,\cdot)\Vert \leq C(\Vert u_{0}\Vert + \Vert d_{n}(\cdot)u_{1}\Vert),\quad t \geq 0.\]
\end{lem}

{\it Proof.}\,Note that the function $v(t,x) := \displaystyle{\int_0^{t}}u(s,x)ds$ is the solution of the problem
\begin{align*}
&v_{tt}-\Delta v +V(x)v=u_1, \quad t>0, \;x\in\Omega,\\
&v(0)=0,\;  v_t(0)=u_0.
\end{align*}
Using the multiplier $m(v) = v_t$ one obtains for $\varepsilon >0$,
\begin{align*}
\Vert v_t\Vert^2 + \Vert \nabla v\Vert^2 + \Vert \sqrt{V(\cdot)}v\Vert^2 &=\Vert u_0 \Vert^2 + 2 \int_{\Omega}u_1(x)v(t,x)dx, \\
&\leq  \Vert u_0 \Vert^2 +  \frac{1}{2\varepsilon} \Vert d_n(\cdot)u_1\Vert^2 + \frac{\varepsilon}{2} \int_{\Omega} \frac{v^2(t,x)}{d_n^2(x)}dx,\quad t>0.
\end{align*}
Applying the Hardy inequality for dimension $n \geq 2$ and choosing a suitable $\varepsilon>0$ the proof of lemma follows due to $u=v_t$.

\hfill
$\Box$

%%%%%%%%%%%%%%%%%%%%%%%%%%%%%%%%%%%%%%%%%%%%%%%%%%%%%%%%section3%%%%%%%%%%%%%%%%%%%%%%%%%%%%%%%%%%%%%%%%%%%%%%%%%%%%%%%%%%%%%%%%%%%%%%%%%%%%%%

\section{Proof of Theorem \ref{t01}}

In this section one proves Theorem \ref{t01} by using Lemmas \ref{10}, \ref{11} and \ref{12}.\\

One first use assumptions {\bf (A-1)} and {\bf (A-2)} and Lemma \ref{10} to get the inequality that
\begin{equation}\label{13}
t E(t) + \frac{n-1}{2}\int_{\Omega}u_{t}(t,x)u(t,x)dx + \int_{\Omega}u_{t}(t,x)(x\cdot\nabla u(t,x))dx \leq J_{0},
\end{equation}
where one has just used the fact that {\bf (A-1)} implies $\sigma\cdot\nu(\sigma) \leq 0$ for each $\sigma \in \partial\Omega$. Thus, it suffices to estimate two ingredients included in \eqref{13} such that
\begin{equation}\label{14}
I_{1}(t) := \left\vert\int_{\Omega}u_{t}(t,x)u(t,x)dx\right\vert,
\end{equation}
\begin{equation}\label{15}
I_{2}(t) := \left\vert\int_{\Omega}u_{t}(t,x)(x\cdot\nabla u(t,x))dx\right\vert.
\end{equation}
$I_{1}(t)$ can be estimated by the use of Lemma \ref{12} soon, and $I_{2}(t)$ can be evaluated by Lemma \ref{11}.\\

{\rm {\bf (I)}.}\, \underline{Let us first get the bound for $I_{1}(t)$.}\\

Indeed, from the Schwarz inequality one has
\[I_{1}(t) \leq \int_{\Omega}\vert u(t,x)\vert\vert u_{t}(t,x)\vert dx\]
\[\leq \left(\int_{\Omega}\vert u_{t}(t,x)\vert^{2}dx \right)^{1/2}\left( \int_{\Omega}\vert u(t,x)\vert^{2}dx \right)^{1/2} = \Vert u_{t}(t,\cdot)\Vert\Vert u(t,\cdot)\Vert.\] 

From \eqref{E1} we see that $\frac{1}{2}\Vert u_{t}(t,\cdot)\Vert^{2} \leq E(t) = E(0)$, so that $\Vert u_{t}(t,\cdot)\Vert \leq \sqrt{2E(0)}$. Thus, by combining Lemma \ref{12} one has
\begin{equation}\label{16}
I_{1}(t) \leq C\sqrt{2E(0)}(\Vert u_{0}\Vert + \Vert d_{n}(\cdot)u_{1}\Vert) \quad t \geq 0.
\end{equation}

{\rm {\bf (II)}.}\,\underline{Let us treat $I_{2}(t)$ to get the decay rate for the local energy.}\\

For this purpose one defines the weight function $\psi(t,x)$ which has the same style introduced in \cite{Ikken}:
\[ \psi(t,x) = \left\{
    \begin{array}{ll}
    \displaystyle{(1 + \vert x\vert -t)}&
    \qquad \vert x\vert \geq t,\quad x \in {\bf R}^{n},\\[0.2cm]
    \displaystyle{(1+t-\vert x\vert)^{-1}}& \qquad \vert x\vert < t,\quad x \in {\bf R}^{n}.
    \end{array} \right.\]
Then, it is easy to check that $\psi \in C^{1}([0,\infty)\times {\bf R}^{n})$ satisfies
\begin{equation}\label{17}
\psi_{t}(t,x) < 0,\quad t > 0,\quad x \in {\bf R}^{n},
\end{equation}
\begin{equation}\label{18}
\psi_{t}(t,x)^{2} - \vert\nabla \psi(t,x)\vert^{2} = 0,\quad t > 0,\quad x \in {\bf R}^{n}.
\end{equation}

Note that \eqref{18} is the so-called Eikonal equation for \eqref{01}. Therefore, it follows from Lemma \ref{11}, $V(x) \geq 0$, \eqref{17} and \eqref{18} that
\[0 \geq \frac{\partial}{\partial t}\left(\psi(t,x)e(t,x)\right) - \nabla\cdot\left(\psi(t,x)u_{t}(t,x)\nabla u(t,x) \right),\,\,t > 0,\,\,x \in \Omega.\]
By integrating both sides above on $[0,t]\times \Omega$ and using the divergence theorem and \eqref{03}, one has the weighted energy estimate such that
\[\int_{\Omega}\psi(t,x)\left(\vert u_{t}(t,x)\vert^{2} + \vert\nabla u(t,x)\vert^{2} + V(x)\vert u(t,x)\vert^{2}\right)dx \]
\begin{equation}\label{19}
\leq \int_{\Omega}(1+\vert x\vert)\left(\vert u_{1}(x)\vert^{2} + \vert\nabla u_{0}(x)\vert^{2} + V(x)\vert u_{0}(x)\vert^{2}\right)dx.
\end{equation} 

Now let us estimate $I_{2}(t)$ basing on \eqref{19}. This is just a modification of \cite[Lemma 2.4]{Ikken}, indeed, let $R > \rho_{0}$ be an arbitrary fixed number.
Set $\Omega_{R} := \Omega \cap B_{R}$. Then, for $t > R$ it follows that
\[I_{2}(t) \leq \int_{\Omega}\vert x\vert\vert u_{t}(t,x)\vert\vert\nabla u(t,x)\vert dx\]
\[\leq R\int_{\Omega_{R}}\vert u_{t}(t,x)\vert\vert\nabla u(t,x)\vert dx + \int_{\vert x\vert \geq R}\vert x\vert\vert u_{t}(t,x)\vert\vert\nabla u(t,x)\vert dx\]
\[\leq \frac{R}{2}\int_{\Omega_{R}}\left(\vert u_{t}(t,x)\vert^{2} + \vert\nabla u(t,x)\vert^{2}\right) dx + \int_{\vert x\vert \geq t}\vert x\vert\vert u_{t}(t,x)\vert\vert\nabla u(t,x)\vert dx + \int_{t \geq \vert x\vert \geq R}\vert x\vert\vert u_{t}(t,x)\vert\vert\nabla u(t,x)\vert dx\]
\[\leq \frac{R}{2}\int_{\Omega_{R}}\left(\vert u_{t}(t,x)\vert^{2} + \vert\nabla u(t,x)\vert^{2} + V(x)\vert u(t,x)\vert^{2}\right) dx\]
\[+ \int_{\vert x\vert \geq t}(\vert x\vert - t)\vert u_{t}(t,x)\vert\vert\nabla u(t,x)\vert dx + t\int_{\vert x\vert \geq t}\vert u_{t}(t,x)\vert\vert\nabla u(t,x)\vert dx\]
\[+ t\int_{t \geq \vert x\vert \geq R}\vert u_{t}(t,x)\vert\vert\nabla u(t,x)\vert dx\]
\[\leq \frac{R}{2}\int_{\Omega_{R}}\left(\vert u_{t}(t,x)\vert^{2} + \vert\nabla u(t,x)\vert^{2} + V(x)\vert u(t,x)\vert^{2}\right) dx + \frac{1}{2}\int_{\vert x\vert \geq t}(1+\vert x\vert - t)\left(\vert u_{t}(t,x)\vert^{2} + \vert\nabla u(t,x)\vert^{2}\right)dx \]
\[+ \frac{t}{2}\int_{\vert x\vert \geq t}\left(\vert u_{t}(t,x)\vert^{2} + \vert\nabla u(t,x)\vert^{2}\right) dx + \frac{t}{2}\int_{t \geq \vert x\vert \geq R}\left(\vert u_{t}(t,x)\vert^{2} + \vert\nabla u(t,x)\vert^{2}\right) dx.\]

Then, one obtain
\[I_{2}(t) \leq \frac{R}{2}\int_{\Omega_{R}}\left(\vert u_{t}(t,x)\vert^{2} + \vert\nabla u(t,x)\vert^{2} + V(x)\vert u(t,x)\vert^{2}\right) dx\]
\[+ \frac{1}{2}\int_{\vert x\vert \geq t}(1+\vert x\vert - t)\left(\vert u_{t}(t,x)\vert^{2} + \vert\nabla u(t,x)\vert^{2}+ V(x)\vert u(t,x)\vert^{2}\right)dx\]
%\[ + \frac{t}{2}\int_{\vert x\vert \geq t}\left(\vert u_{t}(t,x)\vert^{2} + \vert\nabla u(t,x)\vert^{2}+ V(x)\vert u(t,x)\vert^{2}\right) dx\]
%%%
\[+ \frac{t}{2}\int_{\vert x\vert \geq R}\left(\vert u_{t}(t,x)\vert^{2} + \vert\nabla u(t,x)\vert^{2}+ V(x)\vert u(t,x)\vert^{2}\right) dx\]
\[\leq \frac{R}{2}\int_{\Omega_{R}}\left(\vert u_{t}(t,x)\vert^{2} + \vert\nabla u(t,x)\vert^{2} + V(x)\vert u(t,x)\vert^{2}\right) dx \]
%%%
\[+ \int_{\vert x\vert \geq t}\psi(t,x)e(t,x)dx + t\int_{\vert x\vert \geq R}e(t,x)dx ,\]
which implies
\[I_{2}(t) \leq RE_{R}(t) + \int_{\Omega}\psi(t,x)e(t,x)dx + t\int_{\vert x\vert \geq R}e(t,x)dx \]
\begin{equation}\label{20}
\leq RE_{R}(t)  + \int_{\Omega}(1+\vert x\vert)e(0,x)dx + t\int_{\vert x\vert \geq R}e(t,x)dx
\end{equation}
because of \eqref{19}, where
\begin{equation}\label{23}
\int_{\Omega}(1+\vert x\vert)e(0,x)dx = \int_{\Omega}(1+\vert x\vert)\left(\vert u_{1}(x)\vert^{2} + \vert\nabla u_{0}(x)\vert^{2} + V(x)\vert u_{0}(x)\vert^{2}\right)dx =:I_{0}.
\end{equation}

Let us prove Theorem \ref{t01} at a stroke.\\

{\it \underline{Proof of Theorem \ref{t01}.}}\,\,Let $R > \rho_{0}$ be an arbitrary fixed number, and take $t > R$. One first gets the inequality from \eqref{13} that
\begin{equation}\label{21}
tE_{R}(t) + t\int_{\vert x\vert \geq R}e(t,x)dx \leq \frac{n-1}{2}I_{1}(t) + I_{2}(t) + J_{0}.
\end{equation}
Because of \eqref{16}, \eqref{20} and \eqref{21} one can get
\[tE_{R}(t) + t\int_{\vert x\vert \geq R}e(t,x)dx \leq J_{0} + C\frac{n-1}{2}\sqrt{E(0)}(\Vert u_{0}\Vert + \Vert d_{n}(\cdot)u_{1}\Vert)\]
\[+ RE_{R}(t) + t\int_{\vert x\vert \geq R}e(t,x)dx + \frac{1}{2}\int_{\Omega}(1+\vert x\vert)\left(\vert u_{1}(x)\vert^{2} + \vert\nabla u_{0}(x)\vert^{2} + V(x)\vert u_{0}(x)\vert^{2}  \right)dx,\]
which implies the desired decay estimate for the local energy:
\[(t-R)E_{R}(t) \leq J_{0} + C\sqrt{E(0)}(\Vert u_{0}\Vert + \Vert d_{n}(\cdot)u_{1}\Vert)  + \frac{I_0}{2}.\]

Note that the exterior energy $\displaystyle{\int_{\vert x\vert \geq R}}e(t,x)dx$ in the region $\vert x\vert \geq R$ can be cancelled with $t$-times nicely in the computations above. That is, one can not have any information on decay in time to the exterior energy.  
\hfill
$\Box$

\begin{rem}{\rm If one relies on the generalized assumption \eqref{ike-501} in place of (A-2), from the proof above the following additional quantity must be estimated in our method:
\[\frac{\gamma}{2}\int_{0}^{\infty}\int_{\Omega}V(x)\vert u(s,x)\vert^{2}dxds < +\infty.\]  
Such a estimate may be extremely difficult and that is not our goal at this work.}
\end{rem}
%%%%%%%%%%%%%%%%%%%%%%%%%%%%%%%%%%%%%%%%%%%%%%%%%%%%%%%%%%%%%%%%%%%%%%%%%%%%%%%%%%%%%%%%%%%%%%%%%%%%%%%%%%%%%%%%%%%%%%%%%%%%%%

\section{Concluding remark}

Let $t > R > \rho_{0}$. From \eqref{19} and \eqref{23} one has
%\begin{equation}
\[\int_{\vert x\vert \geq t}(1+\vert x\vert -t)e(t,x)dx \leq \int_{\Omega}(1+\vert x\vert)e(0,x)dx = I_{0}.\]
%\end{equation}
Then, for any fixed small $\varepsilon > 0$ one has
\[\int_{\vert x\vert \geq (1+\varepsilon)t}(1+\vert x\vert -t)e(t,x)dx \leq I_{0}.\]
This implies
%\begin{equation}\label{22}
\[\int_{\vert x\vert \geq (1+\varepsilon)t}e(t,x)dx \leq \frac{I_{0}}{1+\varepsilon t}.\]
%\end{equation}

While one knows the energy conservation identity such that $E(t) = E(0)$. Thus, one has a decomposition of the total energy such that
\[\int_{\vert x\vert \geq (1+\varepsilon)t}e(t,x)dx + \int_{R \leq \vert x\vert \leq (1+\varepsilon)t}e(t,x)dx + E_{R}(t) = E(0).\]

Therefore, one can observe the energy concentration integral quantity such that
\begin{equation}\label{22}
\int_{R \leq \vert x\vert \leq (1+\varepsilon)t}e(t,x)dx = E(0) + O(t^{-1}) + O(\frac{1}{1+\varepsilon t}), \quad (t \gg 1)
\end{equation}
by using Theorem 1.2. One observes that \eqref{22} may express a typical wave property from the viewpoint of the energy propagation under the non-compact support assumption on the initial data.

%%%%%%%%%%%%%%%%%%%%%%%%%%%%%%%%%%%%%%%%%%%%%%%%%%%%%%%%%%%%%%%%%%%%%%%%
\par
\vspace{0.5cm}
\noindent{\em Acknowledgement.}
\smallskip
The author deeply would lik to thank my friend Ruy Coimbra Char\~ao (UFSC, Brazil) for his useful comments, suggestions, and careful reading on the first draft. The work of the author was supported in part by Grant-in-Aid for Scientific Research (C) 20K03682  of JSPS. 

%%%%%%%%%%%%%%%%%%ހreferenceހ%%%%%%%%%%%%%%%%%%%%%%%%


\begin{thebibliography}{99}
\bibitem{AGPP} A. Arnold, S. Geevers, I. Perugia and D. Ponomarev, On the exponential time-decay for the one-dimensional wave equation with variable coefficients, Comm. Pure Appl. Anal., 21(10) (2022), 3389--3405. doi:10.3934/cpaa.2022105

\bibitem{EVR}E. Bisognin, V. Bisognin and R. C. Char\~ao, Uniform stabilization for elastic waves system with highly nonlinear localized dissipation, Portugaliae Mathematica 60 (2003), 99--124.

%\bibitem{B} J.-M. Bouclet, N. Burq, Sharp resolvent and time decay estimates for dispersive equations on asymptotically Euclidean backgrounds, arXiv: 1810.01711v1 [math.AP] 3 Oct 2018.

\bibitem{B} N. Burq, D\'ecroissance de l'\'energie locale de L'\'equation des ondes pour le probleme ext\'erieur et absence de r\'esonance au voisinage du r\'eel, Acta Math. 180 (1998), 1-29.

\bibitem{C-Ruy} R. C. Char\~ao, On the principle of limiting amplitude
for perturbed elastic waves in 3D, Bolletino U.M.I. (7)  10-B (1996), 781-797. 

\bibitem{C-Perla} R. C. Char\~ao and G. P. Menzala, On Huygen's principle and perturbed elastic waves, Diff. and  Integral Eqns. 5 (1992), 631-646.

\bibitem{CI} R. C. Char\~ao and R. Ikehata, A note on decay rates of the local energy for wave equations with Lipschitz wavespeeds, J. Math. Anal. Appl. 483 (2020), 123636.

\bibitem{DS} W. Dan and Y. Shibata, On a local energy decay of solutions of a dissipative wave equation, Funkcialaj Ekvacioj 38 (1995), 545--568.

\bibitem{GV} V. Georgiev and N. Visciglia, $L^{\infty}$-$L^{2}$ weighted estimate for the wave equation with potential, Atti Accad. Naz. Lincei Cl. Sci. Fis. Mat. Natur. Rend. Lincei (9) Mat. Appl. 14 (2003), no.2, 109--135.


\bibitem{Ikehat} R. Ikehata, Local energy decay for linear wave equations with non-compactly supported initial data, Math. Meth. Appl. Sci., 27 (2004) 1881-1892.

\bibitem{Ikehat2} R. Ikehata, Local energy decay for linear wave equations with variable coefficients, J. Math. Anal. Appl. 306 (2005) 330-348.

\bibitem{IM} R. Ikehata and T. Matsuyama, $L^{2}$-behaviour of solutions to the linear heat and wave equations in exterior domains, Sci. Math. Japonicae 55 (2022), 33--42.

\bibitem{Ikken} R. Ikehata and K. Nishihara, Local energy decay for  wave equations with initial data decaying slowly near infinity, Gakuto Iternational Series - Mathematical Sciences and Applications,  The $5^{th}$ East Asia PDE Conference, 22 (2005) 265-275.

\bibitem{LP} P. D. Lax and R. S. Phillips, Scattering theory, Revised Edition. Academic Press, New York, 1989.

\bibitem{maraw} C. S. Morawetz, The decay of solutions of the exterior initial-boundary value problem
for the wave equation, Comm. Pure Appl. Math., 14 (1961) 561-568.

\bibitem{morast} C. S. Morawetz, J. V. Ralston and W. A. Strauss, Decay of solutions of the wave equation outside nontrapping obstacles, Comm. Pure Appl. Math., 30 (1977) 447-508.

\bibitem{Mur} L. A. Murav$\breve{e}$, The wave equation in an unbounded domain with a star-shaped boundary, Soviet Math. Dokl. 38 (1989), 527-530.

\bibitem{Nu} R. S. O. Nunes and W. D. Bastos, Energy decay for the linear Klein-Gordon equation and boundary control, J. Math. Anal.
Appl. 414 (2014) 934-944.

\bibitem{R} J. Ralston, Solutions of the wave equation with localized energy, Comm. Pure Appl. Math. 22 (1969), 807-823.

\bibitem{S} J. Shapiro, Local energy decay for Lipschitz wavespeeds, Comm. Partial Diff. Eqns 43, Issue 5 (2018), 839-858. DOI: 10.1080/03605302.2018.1475491.

\bibitem{ST} Y. Shibata and Y. Tsutsumi, On a global existence theorem of small amplitude solutions for nonlinear wave equations in an exterior domain, Math. Z. 191 (1986), 165--199.

\bibitem{T} H. Tamura, On the decay of local energy for wave equations with time-dependent potentials, J. Math. Soc. Japan 33, No. 4 (1981), 605--618.

\bibitem{TY} G. Todorova and B. Yordanov, Critical exponent for a nonlinear wave equation with damping, J. Diff. Eqns 174, (2001), 464--489.

\bibitem{Va} B. R. Vainberg, On the short wave asymptotic behavior of solutions of stationary problems and the asymptotic behavior as $t \to \infty$ of solutions of nonstationary problems, Russian Math. Survey 30 (1975), 1-58.

\bibitem{V} G. Vodev, Local energy decay of solutions to the wave eqution for short-range potentials, Asymptotic Anal. 37 (2004), 175--17.

\bibitem{Zac-2} E. C. Zachmanoglou, The decay of solutions of the initial-boundary value problem for the wave equation in unbounded regions,
Arch. Rational Mech. Anal. 14 (1963) 312-325.

%\bibitem{Zac} E. C. Zachmanoglou, The decay of the initial-boundary value problem for hyperbolic equations,
% J. Math. Anal. Appl. 13 (1966) 504-515.

%\bibitem{V} B. R. Vainberg, On the short wave asymptotic behavior of solutions of stationary problems and the asymptotic behavior as $t \to \infty$ of solutions of nonstationary problems, Russian Math. Survey 30 (1975), 1-58.
\end{thebibliography}
\end{document}